\documentclass{article}

\title{A Generalized Next-Closure Algorithm -- Enumerating Semilattice Elements\\
  from a Generating Set}
\author{
  \small Daniel Borchmann\\
  \small TU Dresden\\
  \small Faculty of Mathematics and Sciences\\
  \small Institute for Algebra\\
  \small \texttt{daniel.borchmann@mailbox.tu-dresden.de}}

\usepackage{mydefs,todonotes}

\begin{document}

\maketitle

\begin{abstract}
  We present a generalization of the well known Next-Closure algorithm working on semilattices. We
  prove the correctness of the algorithm and apply it on the computation of the intents of a formal
  context.
\end{abstract}

\section{Introduction}

Next-Closure is one of the best known algorithms in Formal Concept Analysis~\cite{ganter1999formal}
to compute the concepts of a formal context.  In its general form it is able to efficiently
enumerate the closed sets of a given closure operator on a finite set.  This generality might be a
drawback concerning efficiency compared to other algorithms like
Close-by-One~\cite{conf/cla/VychodilKO10,conf/iccs/Andrews09}, but widens its field of applications.
However, there are still applications where Next-Closure might be useful, but is not applicable,
because an closure operator on a finite set is not explicitly available.  One such example might be
the computation of concepts of a fuzzy formal context~\cite{DBLP:conf/icfca/BelohlavekV06}.  Even
worse, if a closure operator is given, but not on a finite set, Next-Closure is not directly
applicable as well.  In those cases most often an ad hoc variation of Next-Closure can be
constructed.  The aim of this paper is to provide a generalization of Next-Closure which covers
those cases, and even goes beyond them.

As it turns out, Next-Closure is not about enumerating closed sets of a closure operator, even not
on an abstract ordered set.  The algorithm is merely about enumerating elements of a certain
semilattice, given as an operation together with a generating set.  This observation is somewhat
surprising, but, as we shall see, quite natural.

This paper is organized as follows.  First of all we shall revisit the original version of
Next-Closure, together with the basic definitions.  Then we present our generalized version working
on semilattices, together with a complete proof of its correctness.  Then we show how this
generalized form is indeed a generalization of the original Next-Closure.  Additionally, we discuss
a new algorithm for enumerating the intents of a given formal context.  Finally, we give some
outlook on further questions which might be interesting within this line of research.

\section{The Next-Closure Algorithm}

Before we are going to discuss our generalized form of Next-Closure, let us revisit the original
version as it is given in~\cite{ganter1999formal,fca:Ganter:1984}.  To make our discussion a bit
more consistent, we shall make one minor modifications to the presentation given here, which will be
explicitly mentioned.

Let $M$ be a finite set and let $c:\subsets{P}\to\subsets{P}$ be a function such that
\begin{enumerate}[a) ]
\item $c$ is idempotent, i.e.\ $c(c(A)) = c(A)$ for all $A\subseteq M$,
\item $c$ is monotone, i.e.\ if $A\subseteq B$, then $c(A)\subseteq c(B)$ for all $A,B\subseteq M$, and
\item $c$ is extensive, i.e.\ $A\subseteq c(A)$ for all $A\subseteq M$.
\end{enumerate}
A set $A\subseteq M$ is called \emph{closed (with respect to $c$)} if $A = c(A)$, and the
\emph{image of $c$} is defined as
\[ c[\subsets{M}] := \set{c(A)\mid A\subseteq M}. \]

Without loss of generality, let $M = \set{1,\ldots,n}$ for some $n\in\NN$.  For two sets $A,B\in
c[\subsets{M}]$ with $A\neq B$ and $i\in M$ we say that $A$ is \emph{lectically smaller} than $B$
\emph{at position $i$} if and only if
\[ i = \min(A\mathrel{\Delta}B) \quad\text{and}\quad i\in B, \]
where $A\mathrel{\Delta}B = (A\setminus B)\cup(B\setminus A)$ is the symmetric difference of $A$ and
$B$.  We shall write $A \prec_i B$ if $A$ is lectically smaller than $B$ at position $i$.  Finally,
we say that $A$ is \emph{lectically smaller} than $B$, for $A,B\in c[\subsets{M}]$, if $A = B$ or $A
\prec_i B$ for some $i\in M$ and we shall write $A \preceq B$ in this case.

It has to be noted that, in contrast to our definition, the lectic order is normally
defined for all sets $A,B\subseteq M$ in the very same spirit as given above.  However, as
we shall see, this is not necessary, wherefore we have restricted our definition to closed
sets only.

Now let us define for $A\in c[\subsets{M}]$ and $i\in M$
\[ A\oplus i := c(\set{j\in A\mid j < i} \cup \set{i}). \]
Then we have the following result.

\begin{Theorem}[Next-Closure]
  Let $A\in c[\subsets{M}]$.  Then the next closed set $A^+\in c[\subsets{M}]$ after $A$ with
  respect to the lectic order $\preceq$, if it exists, is given by
  \[ A^+ = A\oplus i \]
  with $i\in M$ being maximal with $A <_i A\oplus i$.
\end{Theorem}

This is the original version of Next-Closure, as it is given in~\cite{ganter1999formal}.

Now let us have a closer look on the definition of $\oplus$.  The set $A\oplus i$ can be seen as the
smallest closed set containing both $\set{j \in A\mid j < i}$ and $\set{i}$, or equivalently, both
$c(\set{j \in A\mid j < i})$ and $c(\set{i})$.  This means that we can rewrite $A\oplus i$ as
\[ A\oplus i = c(\set{j \in A\mid j < i}) \vee c(\set{i}), \]
%\
where $X\vee Y$ is the smallest closed sets containing both $X, Y\in c[\subsets{M}]$, the
\emph{supremum of $X$ and $Y$}, which is simply given by $X \vee Y = c(X\cup Y)$.  This observation
allows us to consider Next-Closure on abstract algebraic structures with a binary operation $\vee$
with some certain properties.  To do so we need a more general notion of $c(\set{i})$, since we do
not necessarily deal with subsets, and a more general notion of $\set{j\in A\mid j<i}$, which
likewise might not be expressible in a more general setting.  Finally, we need to find a starting
point for our enumeration, which is $c(\emptyset)$ in the original description of Next-Closure, but
may vary in other cases.  Luckily, all this is possible and quite natural, as we shall see in the
next section.

\section{Generalizing Next-Closure for Semilattices}

The aim of this section is to present a generalization of the Next-Closure algorithm that works on
semilattices. For this recall that a semilattice $\underline{L} = (L,\vee)$ is an algebraic
structure with a binary operation $\vee$ which is associative, commutative and idempotent. It is
well known that with
\[ x \le_{\underline{L}} y \,:\!\iff x\vee y = y, \]
with $x,y\in L$.  An order relation on $L$ is defined such that for every two elements
$x,y$ the element $x\vee y$ is the least upper bound of both $x$ and $y$ with respect to
$\le_{\underline{L}}$.

For the remainder of this section let $\underline{L} = (L,\vee)$ be an arbitrary but fixed
semilattice. Furthermore, let $(x_i\mid i\in I)$ be an enumeration of a finite generating set
$\set{x_i\mid i\in I}\subseteq L$ of $\underline{L}$. Finally, let $\le_I$ be a total order on $I$.

\newcommand{\leL}{\le_{\underline{L}}}

\begin{Definition}
  Let $a,b\in L$ and let $i\in I$. Set
  \[
  \Delta_{a,b} := \set{j\in I
    \mid
    (x_j\le_{\underline{L}} a \text{ and } x_j\not\le_{\underline{L}} b)
    \text{ or }
    (x_j \not\le_{\underline{L}} a \text{ and } x_j\le_{\underline{L}} b)}.
  \]
  We then define
  \[
    a <_i b \,:\!\iff i = \min\Delta_{a,b} \text{ and } x_i\le_{\underline{L}} b.
  \]
  Furthermore we write $a < b$ if $a <_i b$ for some $i\in I$ and write $a\le b$ if $a = b$ or $a <
  b$.
\end{Definition}

One can see the similarity of this definition to the one of the lectic order.  Here, the set
$\Delta_{a,b}$ generalizes $a\mathrel{\Delta}b$ and $x_i\leL b$ somehow represents the fact that
$i\in b$, or equivalently $\set{i}\subseteq b$, in the special case of $L = \subsets{M}$ and $i\in
M$.

Note that if $a <_i b$ and $k\in I$ with $k<_I i$, then
\[ x_k\leL a \iff x_k\leL b. \]
This observation is quite useful and will be used in some of the proofs later on.

The first thing we want to consider now are two easy results stating that $\le$ is a total order
relation on $L$ extending $\leL$.

\begin{Lemma}\label{lem:lectic is order}
  The relation $<$ is irreflexive and transitive.  Furthermore, for every two elements $a,b\in L$
  with $a\neq b$, it is either $a < b$ or $b < a$.

  \begin{Proof}
    If $a = b$, then the set $\Delta_{a,b}$ defined above is empty, therefore we cannot have $a <_i
    a$ for some $i\in I$.  This shows the irreflexivity of $<$.  Let us now consider the
    transitivity of $<$.  For this let $a,b,c\in L$, $i,j\in I$ and suppose that $a<_ib$ and $b<_j
    c$.  We have to show that $a < c$.  Let us consider the following cases.

    \textit{Case $i<_I j$.} We have $x_i\not\leL a$ and $x_i\leL b$ because of $a<_i b$. Due to $i
    <_I j$ it follows that $x_i\leL c$.  Suppose that there exists $k\in I$, $k <_I i$ with $x_k
    \leL a$ and $x_k\not\leL c$.  Then if $x_k\not\leL b$ we would have $x_k\not\leL a$ because $k
    <_I i$, a contradiction.  But if $x_k\leL b$, then $x_k\leL c$ because of $k <_I i <_I j$, again
    a contradiction.  Thus we have shown that $a < c$.

    \textit{Case $j <_I i$.} We have $x_j\not\leL b$, $x_j\leL c$ because of $b <_j c$.  Due to $j
    <_I i$ it follows that $x_j \not\leL a$.  Now if there were a $k\in I$, $k <_I j$ with $x_k\leL
    a$ and $x_k\not\leL c$, then $x_k\leL b$ would imply $x_k\leL c$ and $x_k\not\leL b$ would imply
    $x_k\not\leL a$, analogously to the first case, a contradiction.  Hence such a $k$ cannot exist
    and $a < c$.

    \textit{Case $i = j$.} This cannot occur since otherwise $x_i\leL b$, because of $a<_i b$, and
    $x_i\not\leL b$, because of $b<_i c$, a contradiction.

    Overall we have shown that $a < c$ in any case and therefore $<$ is a transitive relation.

    Finally let $a,b\in L$ with $a\neq b$.  Then because $\set{x_i\mid i\in I}$ is a generating set,
    the set $\Delta_{a,b}$ is not empty, since otherwise $a = b$, and with $i := \min\Delta_{a,b}$
    we either have $a <_i b$ if $x_i\leL b$ and $b <_i a$ otherwise.
  \end{Proof}
\end{Lemma}

\begin{Lemma}\label{lem:lectic extends original order}
  Let $a,b\in L$ with $a\le_{\underline{L}} b$. Then $a\le b$. In particular, if $a \leL c$ and
  $b\leL c$ for $a,b,c\in L$, then $a\vee b \le c$.

  \begin{Proof}
    We show $x_i\leL a\implies x_i\leL b$ for all $i\in I$.  This shows $b\not<a$, hence $a\le b$ by
    Lemma~\ref{lem:lectic is order}.  Now if $x_i\leL a$, then because of $a\leL b$ we see that
    $x_i\leL b$ and the claim is proven.
  \end{Proof}
\end{Lemma}

The next step towards a general notion of Next-Closure is to provide a generalization of $\oplus$.

\begin{Definition}
  Let $a\in L$ and $i\in I$. Then define
  \[ a \oplus i := \bigvee_{\substack{j <_I i\\x_j\le_{\underline{L}} a}} x_j \vee x_i. \]
\end{Definition}

With all these definitions at hand we are now ready to prove the promised generalization.
For this, we generalize the proof of Next-Closure as it is given in~\cite[page
67]{ganter1999formal}.

\begin{Lemma}\label{lem:helper}
  Let $a,b\in L$ and $i,j\in I$. Then the following statements hold:
  \begin{enumerate}[i) ]
  \item\label{lem:helper:1} $a <_i b, a <_j c, i <_I j \implies c <_i b$.
  \item\label{lem:helper:2} $a < a \oplus i$ if $x_i\not\le_{\underline{L}} a$.
  \item\label{lem:helper:3} $a <_i b \implies a\oplus i \le b$.
  \item\label{lem:helper:4} $a <_i b \implies a <_i a\oplus i$.
  \end{enumerate}

  \begin{Proof}
    \begin{enumerate}[i) ]
    \item It is $x_i\not\leL a$ and due to $i<_I j$ we get $x_i\not\leL c$ as well.  Furthermore,
      $x_i\leL b$ because of $a <_i b$.  Now if there would exist a $k\in I$ with $k <_I i$ such
      that $x_k\leL c, x_k\not\leL b$, then $x_k\leL a$ because of $k <_I i$ and $x_k\not\leL a$
      because of $k <_I i <_I j$, a contradiction. With the same argumentation a contradiction
      follows from the assumption that there exists a $k\in I$, $k<_I i$ with $x_k\not\leL c,
      x_k\leL b$.  In sum we have shown $c <_i b$, as required.
    \item We have $x_i\not\leL a$ and $x_i\leL a\oplus i$.  Furthermore, for $k\in I$, $k <_I i$ and
      $x_k\leL a$ we have $x_k\leL a\oplus i$ by definition.  This shows $a < a\oplus i$.
    \item Let $a <_k b$ for some $k\in I$. Then $\bigvee_{j <_I k,x_j\leL a}x_j \leL b$ and $x_k\leL
      b$, hence with Lemma~\ref{lem:lectic extends original order} we get $a\oplus k\le b$.
    \item Let $a <_i b$.  Then $x_i\not\leL a$ and with~(\ref{lem:helper:2}) we get $a < a\oplus
      i$. By~(\ref{lem:helper:3}), $a\oplus i \le b$. If for $k\in I$, $k <_I i$ it holds that
      $x_k\leL a\oplus i$ and $x_k\not\leL a$, then we also have $x_k\leL a\oplus i \leL b$, i.e.\
      $x_k\leL b$, contradicting the minimality of $i$.
    \end{enumerate}
  \end{Proof}
\end{Lemma}

\begin{Theorem}[Next-Closure for Semilattices]\label{thm:semiNC}
  Let $a\in L$. Then the next element $a^+\in L$ with respect to $<$, if it exists, is given by
  \[ a^+ = a\oplus i \]
  with $i\in I$ being maximal with $a <_i a\oplus i$.

  \begin{Proof}
    Let $a^+$ be the next element after $a$ with respect to $<$. Then $a <_i a^+$ for some
    $i\in I$ and by Lemma~\ref{lem:helper}.\ref{lem:helper:4} we get $a <_i a\oplus i$ and
    with Lemma~\ref{lem:helper}.\ref{lem:helper:3} we see $a\oplus i \le a^+$, hence
    $a\oplus i = a^+$. The maximality of $i$ follows from
    Lemma~\ref{lem:helper}.\ref{lem:helper:1}.
  \end{Proof}
\end{Theorem}

To find the correct element $i\in I$ such that $a^+ = a\oplus i$ can be optimized with
Lemma~\ref{lem:helper}.\ref{lem:helper:2}.  Because of this result, only elements $i\in I$ with
$x_i\not\leL a$ have to be considered, a technique which is also known for the original form of
Next-Closure.

However, to make the above theorem practical for enumerating the elements of a certain semilattice,
one has to start with some element, preferably the smallest element in $L$ with respect to $\le$.
This element must also be minimal in $L$ with respect to $\leL$, by Lemma~\ref{lem:lectic extends
  original order}.  Since $\set{x_i\mid i\in I}$ is a generating set of $\underline{L}$, and $a \le
a \vee b$ for all $a,b\in L$, the minimal elements of $L$ with respect to $\leL$ must be among the
elements $x_i$, $i\in I$.  So to find the first element of $\underline{L}$ with respect to $\le$,
find all minimal elements in $\set{x_i\mid i\in I}$ and choose the smallest element with respect to
$\le$ from them.  But because of $x_i \le x_j$ if and only if $j <_I i$, one just has to take the
largest index $j$ of all minimal elements among the $x_i$ to find the smallest element in $L$ with
respect to $\le$.

As a final remark for this section note that the set $\set{x_i\mid i\in I}$ must always include the
$\vee$-irreducible elements of $\underline{L}$.  These are all those elements $a\in L$ that cannot
be represented as a join of other elements, or, equivalently,
\[
  \set{b \in L\mid b <_{\underline{L}} a} = \emptyset
  \quad\text{or}\quad
  \bigvee_{b <_{\underline{L}} a} b <_{\underline{L}} a.
\]
It is also easy to see that the $\vee$-irreducible elements of $\underline{L}$ are also sufficient,
i.e.\ they are a generating set of $\underline{L}$.

%\todo[inline]{Investigate case when $a\oplus i$ with maximal possible $i$ is always the next
%  element}

\section{Computing the Intents of a Formal Context}

We have seen an algorithm that is able to enumerate the elements of a semilattice from a given
generating set.  We have also claimed that this is a generalization of Next-Closure, which we want
to discuss in this section.  Furthermore, we want to give another example of an application of this
algorithm, namely the computation of the intents of a given formal context.

Firstly, let us reconstruct the original Next-Closure algorithm from Theorem~\ref{thm:semiNC} and
the corresponding definitions.  For this let $M$ be a finite set and let $c$ be a closure operator
on $M = \set{0,\ldots,n-1}$, say.  We then apply Theorem~\ref{thm:semiNC} to the semilattice
$\underline{P} = (c[\subsets{M}], \vee)$.  We immediately see that ${\le_{\underline{P}}} =
{\subseteq}$ and that $<_i$ is the usual lectic order on $\underline{P}$.  Then the set
\[ \set{c(\set{i}) \mid i\in M}\cup\set{c(\emptyset)} \]
is a finite generating set of $\underline{P}$ and we can define $x_i := c(\set{i})$ and $x_n :=
c(\emptyset)$, i.e.\ $I = \set{0,\ldots,n}$.  For a closed set $A\subseteq M$ and $i\in I$ then
follows
\begin{align*}
  A\oplus i
  &= \bigvee_{\substack{j < i\\x_j\subseteq A}}x_i \vee x_i\\
  &= c(\bigcup_{\substack{j < i\\x_j\subseteq A}} x_j) \vee x_i\\
  &= c(\bigcup\set{c(\set{j})\mid j< i, j\in A}) \vee c(\set{i})\\
  &= c(\set{j\mid j < i, j\in A}) \vee c(\set{i})\\
  &= c(\set{j\mid j < i, j\in A} \cup \set{i})
\end{align*}
which is the original definition of $\oplus$ for Next-Closure. Furthermore, it is $A\oplus n = A$
since $c(\emptyset)\subseteq A$ for each closed set $A$.  We therefore do not need to consider $x_n$
when looking for the next closed set, and indeed, the only reason why $x_n = c(\emptyset)$ has been
included is that it is the smallest closed set in $\underline{P}$.  All in all, we see that
Next-Closure is a special case of Theorem~\ref{thm:semiNC}.

However, for a closure operator $c$ on a finite set $M$ it seems more natural to consider the
semilattice $\underline{P} = (c[\subsets{M}],\cap)$, because the intersection of two closed sets of
$c$ again yields a closed set of $c$.  One sees that ${\le_{\underline{P}}} = {\supseteq}$.  As a
generating set we take the set of $\cap$-irreducible elements $\set{X_i\mid i\in G}$ for some index
set $G$.  Let $A,B\in c[\subsets{M}]$ and let $<_G$ be a linear ordering on $G$.  Then $A < B$ if
and only if there exists $i\in G$ such that
\[
  i = \min\set{j\in G\mid (X_i\supseteq A, X_i \not\supseteq B) \text{ or } (X_i \not\supseteq A,
    X_i\supseteq B)}
  \quad\text{and}\quad
  X_i \supseteq B
\]
and $\oplus$ is just given by
\[ A\oplus i = \bigcap_{\substack{j <_G i\\X_j\supseteq A}}X_j\cap X_i. \]
Now note that $\oplus$ does not need the closure operator $c$ anymore.  This means that if the
computation of $c$ is very costly and the $\cap$-irreducible elements (or a superset thereof) is
known, this approach might be much more efficient.  In general, however, it is not known how to
efficiently determine the $\cap$-irreducible closed sets of $c$.  But if $c$ is given as the
$\cdot''$ operator of a formal context, these irreducible elements can be determined
quickly~\cite{ganter1999formal}.

Let $G$ and $M$ be two finite sets and let $J\subseteq G\times M$.  We then call the triple
$\mathbb{K} := (G,M,J)$ a \emph{formal context}, $G$ the \emph{objects} of the formal context and
$M$ the \emph{attributes} of the formal context.  For $g\in G$ and $m\in M$ we write $g\mathrel{J}m$
for $(g,m)\in J$ and say that object $g$ \emph{has} attribute $m$.

Let $A\subseteq G$ and $B\subseteq M$.  We then define the \emph{derivations} of $A$ and $B$ to be
\begin{align*}
  A' &:= \set{m\in M\mid \forall g\in A: g\mathrel{J}m}\\
  B' &:= \set{g\in G\mid \forall m\in B: g\mathrel{J}m}.
\end{align*}
Then the $\cdot''$ operator is just the twofold derivation of a given set of attributes.  It turns
out that this is indeed a closure operator, and that every closure operator can be represented as a
$\cdot''$ operator of a suitable formal context~\cite{ganter1999formal,Borchmann11}.  The closed
sets of $\cdot''$, i.e.\ all sets $B\subseteq M$ with $B = B''$, are called the \emph{intents} of
$\mathbb{K}$ and shall be denoted by $\Int(\mathbb{K})$.  It is clear from the previous remarks that
$(\Int(\mathbb{K}),\cap)$ is a semilattice.

The advantage of representing a closure operator is that the $\cap$-irreducible elements of
$(\Int(\mathbb{K}),\cap)$ can be directly read off from the format context.  As discussed
in~\cite{ganter1999formal}, the set
\[ \set{\set{g}'\mid g\in G} \]
contains the irreducible elements we are looking for, except $M$.  Furthermore, it is possible to
omit certain objects $g$ from $\mathbb{K}$ without changing $\Int(\mathbb{K})$.  Every object $g\in
G$ can be omitted from $\mathbb{K}$ for which the set $\set{g}'$ is either equal to $M$ or can be
represented as a proper intersection of other sets $\set{g_1}',\ldots,\set{g_n}'$ for some elements
$g_1,\ldots,g_n\in G$.  It is also clear that if there exist two distinct objects $g_1$ and $g_2$
with $\set{g_1}' = \set{g_2}'$, that we can remove one of them without changing $\Int(\mathbb{K})$.
A formal context for which no such objects exist is called \emph{object clarified} and \emph{object
  reduced}.  If $\con{K} = (G,M,J)$ is an object clarified and object reduced formal context, then
the set $\set{\set{g}'\mid g\in G}$ is exactly the set of $\cap$-irreducible intents of $\con{K}$,
except for the set $M$.

The above described algorithm now takes the following form when applied to $(\Int(\con{K}),\cap)$.
As index set we choose the set $G$ of object of the given formal context, ordered by $<_G$.  For
every object $g\in G$ we set $x_g := \set{g}'$.  Then the set $\set{x_g\mid g\in G}$ is a generating
set of the semilattice $(\Int(\con{K})\setminus\set{M},\cap)$, which we want to enumerate (since we
get the set $M$ for free).  For $A$ begin an intent of $\con{K}$ and $g\in G$ we have
\begin{align*}
  A\oplus g
  &:= \bigcap_{\substack{h\in G\\h<_G g}}\set{h}' \cap \set{g}'
\end{align*}
and as the first intent we take $M$.  For an intent $A\subseteq M$ of $\con{K}$ we then have to find
the maximal object (with respect to $<_G$) $g\in G$ such that $A <_g A\oplus g$.  This is equivalent
to $g$ being maximal with $\set{g}'\not\supseteq A$ and $\forall h\in G,h<_G g: \set{h}'\supseteq A
\iff \set{h}'\supseteq A\oplus g$.  However, the direction ``$\Longrightarrow$'' is clear, hence we
only have to ensure
\[ \forall h\in G, h<_Gg: \set{h}'\not\supseteq A \implies \set{h}'\not\supseteq A\oplus g. \]
All these considerations yield the algorithm shown in Listing~\ref{alg:next closure with cap}.  Of
course, the derivations of the form $\set{g}'$ should not be computed every time they are needed but
rather stored somewhere for reuse.

\begin{lstlisting}[float,caption={Compute Next Intent of a Formal Context},label={alg:next closure
    with cap},language=PseudoCode,numbers=none]
define next-intent($\con{K} = (G,M,J)$,$A$)
  for $g\in G$, descending
    if $\set{g}'\not\supseteq A$ then
      let ($B$ := $A\oplus g$)
        if $\forall h\in G, h<_G g, \set{h}'\not\supseteq A: \set{h}'\not\supseteq B$ then
          return $B$
        end if
      end let
    end if
  end for
  return nil
end
\end{lstlisting}

In contrast to the original version of Next-Closure used to compute the intents of a given formal
context, our version only traverses the formal context $\con{K}$ once, when computing $\set{g}'$ for
every $g\in G$.  This might be an advantage if accessing the context is computationally expensive,
but also be a disadvantage if there are much more objects than attributes (which occurs quite often
in practice).  However, in this case one might better compute the extents of the formal context,
i.e.\ all sets $A\subseteq G$ with $A = A''$, which also form a closure system.  Then the algorithm
in Listing~\ref{alg:next closure with cap} would run through the set of attributes to compute the
extent.

Finally, let us consider the time complexity of the new algorithm to compute the intent after $A$
for a given formal context $\con{K} = (G,M,J)$.  If we assume the operations $\cap$ and $\supseteq$
to be constant (i.e.\ independent of the size of both $G$ and $M$, a very optimistic assumption),
then the algorithm from Listing~\ref{alg:next closure with cap} roughly needs $|G|\times(|G|+|G|)$
steps to compute the next intent, i.e.\ has worst time complexity $\mathcal{O}(|G|^2)$.  On the
other hand, if one assumes a naive implementation of $\cap$ and $\supseteq$ taking time
$\mathcal{O}(|M|)$, one has time complexity $|G|\times(|M|+|G|\times|M|+|G|\times|M|)$, hence worst
case complexity $\mathcal{O}(|G|^2\times|M|)$, which is the time complexity of Next-Closure when
used to compute the extents of $\con{K}$.

\section{Conclusion}

We have seen a natural generalization of the Next-Closure algorithm to enumerate elements of a
semilattice from a generating set.  We have proven the algorithm to be correct and applied it to the
standard task of computing the intents of a given formal context, yielding a new algorithm to
accomplish this.  However, there are still some interesting ideas one might want to look at.

Firstly, a variation of the original Next-Closure algorithm is able to compute the stem base of a
formal context, a very compact representation of its implicational knowledge.  It would be
interesting to know whether the generalization given in this paper gives more insight into the
computation, and therefore into the nature of the stem base.  This is, however, quite a vague idea.

Secondly, the new algorithm discussed above to compute the intents of a formal context enumerates
them in a certain order, which might or might not be a lectic one.  Understanding this order
relation might be fruitful, especially with respect to the complexity results obtained recently,
which consider enumerating pseudo-intents of a formal context in lectic
order~\cite{DBLP:conf/icfca/Distel10}.

Finally, and more technically, it might be interesting to look for other applications of our general
algorithm.  As already mentioned in the introduction, the enumeration of fuzzy concepts of a fuzzy
formal context might be worth investigating (and it might be interesting comparing it
to~\cite{DBLP:journals/tfs/BelohlavekBOV10}), but there might be other applications.

\bibliographystyle{plain}
\bibliography{fca}

\end{document}